\documentclass[draft,10pt,oneside,letterpaper]{article}
\usepackage{fullpage}
\usepackage{amsmath}
\usepackage{amssymb}
\usepackage{amsthm}
\usepackage{calrsfs}
\usepackage[dvips,curve]{xypic}
\usepackage{doublespace}

\newcommand{\End}{\operatorname{End}}

\newcommand{\ind}{\operatorname{ind}}
\newcommand{\res}{\operatorname{res}}
\newcommand{\Der}{\operatorname{Der}}
\newcommand{\Hom}{\operatorname{Hom}}
\newcommand{\Ext}{\operatorname{Ext}}

\newcommand{\isom}{\cong}
\newcommand{\tensor}{\otimes}
\newcommand{\gl}{\mathrm{gl}}

\newcommand{\Liealg}[1]{\mathrm{#1}}
\renewcommand{\sl}{\Liealg{sl}}
\newcommand{\osp}{\Liealg{osp}}
\newcommand{\CC}{\mathbb{C}}

\newcommand{\ZZ}{\mathbb{Z}}
\newcommand{\g}{\mathfrak{g}}
\newcommand{\h}{\mathfrak{h}}
\newcommand{\p}{\mathfrak{p}}

\newcommand{\W}{W}
\newcommand{\borel}{\mathfrak{b}}
\newcommand{\n}{\mathfrak{n}}

\newcommand{\veps}{\varepsilon}

\newcommand{\M}{\mathcal{M}}
\newcommand{\Irr}{\operatorname{Irr}}

\newcommand{\Ob}{\operatorname{Ob}}

\newcommand{\U}{U}
\newcommand{\set}[2]{\left\{\,#1\mathrel{\left|\vphantom{#1}\right.}#2\,\right\}}

\newcommand{\dxi}[1]{\frac{\partial}{\partial\xi_#1}}
\newcommand{\std}{\operatorname{std}}

\newcommand{\myto}{\rightsquigarrow}

\newcommand{\Ann}{\operatorname{Ann}}
\newcommand{\height}[1]{\operatorname{ht}(#1)}

\newcommand{\ad}{\operatorname{ad}}

\theoremstyle{plain}
\newtheorem{thm}{Theorem}[section]
\newtheorem{lemma}[thm]{Lemma}
\newtheorem{conjecture}[thm]{Conjecture}
\newtheorem*{claim}{Claim}
\newtheorem{prop}[thm]{Proposition}

\newtheorem{cor}[thm]{Corollary}
\theoremstyle{remark}
\newtheorem*{definition}{Definition}

\title{Blocks of Lie Superalgebras of Type $\W(n)$}
\author{Noam Shomron\thanks{The author was partially supported by TMR Grant ERB--FMRX--CT97--0100 during the summer of~2000.}}

\begin{document}

\maketitle

\addtocounter{section}{-1}
\section{Introduction}

Let $\g$ be a simple, finite-dimensional Lie superalgebra over~$\CC$.
These have been classified by V.~Kac.  Unless $\g$ is a Lie
algebra or a Lie superalgebra of type~$\osp(1,2n)$,
the category of finite-dimensional representations of~$\g$ is not semisimple;
q.v.~\cite{Scheunert}.
This leads to a classification problem.
For example, in~\cite{Jer}, the representation theory of~$\sl(m,n)$ is 
worked out by showing it is wild when $m,n\ge2$, and by giving an
explicit description of the indecomposable finite-dimensional
representations of~$\sl(1,n)$.

When $\g$ is of type~$\W(0,n)$, the irreducible finite-dimensional
$\g$-modules are classified in~\cite{BL}; in this paper, we
investigate finite-dimensional indecomposable modules.
We show that the category of finite-dimensional representations
of~$\g$ is wild (i.e., as hard as classifying pairs of matrices; q.v.\ \S\ref{sec:quivers}) when $\g$ is of type~$\W(0,n)$ with $n\ge3$.  More precisely, the category of finite-dimensional
representations decomposes into blocks parametrised by
$\bigl(\CC/\ZZ\bigr)\times\ZZ_2$, and we show that each block is of wild type.  This is
done by explicitly exhibiting enough extensions between simple
modules.

Secondly, we find the decomposition of the category of
finite-dimensional representations into blocks.  As an application,
using an idea of Maria~Gorelik, we prove that the centre of the
universal enveloping algebra of~$\g$ is trivial.

When $n=2$,
there is a special isomorphism $\W(0,2)\isom\sl(1,2)$, in
which case the representation theory is not wild, and the
indecomposable representations are fully described in~\cite{Jer}.

The results in this paper are also related to results of Nakano~\cite{Nakano}
in the finite-characteristic case, for which he shows that the
restricted universal enveloping algebra has a single block, and
determines the structure of projective modules.

\section{Conventions}

By~$\g$ we shall denote a simple, finite-dimensional Lie superalgebra
over the complex numbers.
Our primary object of study will be the category of
(graded) finite-dimensional $\g$-modules, with even intertwiners.

The universal enveloping algebra $U(\g)$ of~$\g$ is an associative
superalgebra; it satisfies the graded version of the usual universal
property, thus a $\g$-module is the same thing as a left
$U(\g)$-module.

\section{Quivers and representation type}\label{sec:quivers}

A \emph{quiver} is a directed graph, which consists of a set of vertices
connected by various arrows (possibly including multiple arrows
between two vertices, loops, etc.)  Let $A$ be a unital $\CC$-algebra, and denote
by~$\M$ some category of modules over~$A$.  Denote by~$\Irr\M$ the
set of isomorphism classes of irreducible objects.  The
\emph{Ext-quiver} of~$\M$ is defined to be the quiver whose set of
vertices is $\Irr\M$ and where the number of arrows from $[S_1]$
to~$[S_2]$ is equal to~$\dim\Ext^1(S_1,S_2)$.  This is a combinatorial 
invariant of~$\M$, whose structure gives information about the
representation type.  In particular, we have the following theorem,
proven in~\cite{Virasoro}:

\begin{prop}\label{wildness}
Let $A$ be an algebra, and let $Q$ be a subquiver of the
$\Ext$-quiver of~$A$.
If $Q$ is a connected quiver containing no path of length~$2$, then
there exists a fully faithful functor from the category of
representations of~$Q$ to the category of
$A$-modules.  In particular, the set of isomorphism classes of
indecomposable representations of~$Q$ embeds into the corresponding
set for~$A$.
\end{prop}

The representation theory of quivers is well-established
(see~\cite{Gabriel} for a comprehensive overview).
In
particular, if the underlying graph of a quiver is not of Dynkin or of 
affine type, then the representation theory of the quiver is
\emph{wild}.
More precisely, a small $\CC$-linear Abelian category~$\M$ is defined to be wild if
there exists a full exact embedding from the category of
finite-dimensional representations of~$\CC\langle x,y\rangle$,
the free associative algebra on two generators, into~$\M$.
This has the consequence that the objects of~$\M$ are unclassifiable
in any finite sense.  For example, if $\M$ is wild, then it is
possible to obtain \emph{any}
finite-dimensional algebra as the endomorphisms of some object.

\section{Definition of $\W(n)$}\label{sec:definition_of_W}

Assume that $n\ge2$.
The finite-dimensional Lie superalgebras $\W(n)\mathrel{:=}\W(0,n)$ may be described as
follows: let $\bigwedge[\xi]=\bigwedge[\xi_1,\ldots,\xi_n]$ be the
Grassmann algebra on $n$ generators.  It is a $2^n$-dimensional
associative $\ZZ$-graded algebra.  We set $\W(n)=\Der\bigwedge[\xi]$.
It is a simple Lie superalgebra of dimension~$n2^n$.
The $\ZZ$-grading on~$\W(n)$ is $\W(n)=\bigoplus_{k=-1}^{n-1}\W_k$.

The component~$\W_0$ is canonically isomorphic to~$\gl(n)$; let us
describe the structure of~$\W(n)$ as a $\gl(n)$-module.  Let $\std$ be 
the standard representation of~$\gl(n)$; then there is an isomorphism
$\W(n) \isom \bigwedge(\std) \tensor \std^*$.  If $\std$ is regarded
as a $\ZZ$-graded vector space lying in degree~$1$, then this is an
isomorphism of $\ZZ$-graded vector spaces, so that $\W_k \isom
\bigwedge^{k+1}(\std) \tensor \std^*$; note that this has at most two
irreducible components.

There is a bijection between simple finite-dimensional $\W(n)$-modules 
and simple finite-dimensional $\gl(n)$-modules, realized as follows:
to an irreducible finite-dimensional $\W(n)$-module $V$, associate the 
$\gl(n)$-module $V^{\W_{\ge1}}$, and, conversely, given a simple
finite-dimensional $\gl(n)$-module $M$, take the irreducible quotient
of the induced module $U(\W(n))\tensor_{U(\W_{\ge0})}M$ (q.v.~\cite{BL}).

It is well-known that the irreducible finite-dimensional
representations of~$\gl(n)$, in the ungraded case, are parametrised by 
weights $\lambda=(\lambda_1,\ldots,\lambda_n)$, where
$\lambda_i-\lambda_j\in\ZZ_{\ge0}$ for $i\le j$.  In the super case,
this is true up to parity reversal.
We fix, once and for all, the Cartan subalgebra~$\h$ consisting of diagonal matrices
in~$\gl(n)$; it is also a Cartan subalgebra of~$\W(n)$.  Weights are
written down with respect to the basis $\{\veps_1,\ldots,\veps_n\}$
of~$\h^*$ dual to~$\{E_{11},\ldots,E_{nn}\}$.

\section{Kac modules}

The Lie superalgebra~$\W(n)$ is defined in
Section~\ref{sec:definition_of_W}.  It contains subalgebras
$\gl(n) \subset \sl(1,n) \subseteq \W(n)$, where we will fix the
identifications $\gl(n)\isom\W_0$ and
$\sl(1,n)\isom\W_{-1}\oplus\W_0\oplus\operatorname{span}\{\,\xi_i E
\mid i=1,\ldots,n\,\}$, where $E=\sum_{i=1}^n\xi_i\dxi{i}$ is the
Euler vector field.  Sometimes we will omit the parameters and write
$\W$ for $\W(n)$,\ $\sl$ for $\sl(1,n)$, and $\gl$ for $\gl(n)$.
Let $\g$ be either $\W(n)$ or $\sl(1,n)$.
We fix the Cartan
subalgebra~$\h=\h_0$ consisting of diagonal matrices in~$\gl(n)$, and
consider all weights with respect to~$\h$.  Weights are written down
with respect to the basis $\{\veps_1,\ldots,\veps_n\}$ of~$\h^*$ dual
to~$\{E_{11},\ldots,E_{nn}\}$; we then have $\h^*\isom\CC^n$ and
the root lattice $Q\isom\ZZ^n$.  We fix the Borel subalgebra~$\borel_0\subset\gl(n)$ 
of upper-triangular matrices, which has positive roots
$\{\,\veps_i-\veps_j \mid 1\le i<j\le n\,\}$.  The corresponding set
of highest weights of finite-dimensional irreducible representations
(modulo parity if we consider graded representations)
is
$$\Lambda^+ = \set{(\lambda_1,\ldots,\lambda_n)}{\forall1\le i,j\le
n\;\lambda_i-\lambda_j\in\ZZ_{\ge0}}.$$
Define the Borel subalgebra $\borel_1$ of~$\g$ by
$\borel_1=\borel_0\oplus\g_{\ge1}$.

Let $L_\lambda$ denote the simple, finite-dimensional $\gl(n)$-module
with highest weight~$\lambda$.
\begin{definition}
The \emph{Kac module} $K(\lambda)=K_\lambda$ is the induced representation
\[
K_\lambda = \ind_{\g{\ge0}}^\g L_\lambda = U(\g)\tensor_{U(\g_{\ge0})}L_\lambda
\]
of~$\g$, where $\g_{\ge1}$ acts trivially on~$L_\lambda$.
\end{definition}
The module~$K(\lambda)$ is finite-dimensional, indecomposable, and has highest
weight~$\lambda$ with respect to the Borel subalgebra $\borel_1$ of~$\g$.

There is a bijection between simple finite-dimensional $\g$-modules 
and simple finite-dimensional $\gl(n)$-modules, realised as follows:
to an irreducible finite-dimensional $\g$-module $V$, associate the 
$\gl(n)$-module $V^{\g_{\ge1}}$, and, conversely, given a simple
finite-dimensional $\gl(n)$-module $M$, extend it to~$\g_{\ge0}$, then 
take the unique irreducible quotient
of the induced module $U(\g)\tensor_{U(\g_{\ge0})}M$.  Denote the
unique irreducible quotient of $K(\lambda)$ by~$S(\lambda)$.

For $\g=\W(n)$, the irreducible finite-dimensional representations are
determined explicitly in~\cite{BL}.   
For a generic weight~$\lambda$, the representation~$K_\lambda$ is
irreducible.

For $\sl(1,n)$, the situation is similar, and we will need the
following facts (see~\cite{Kac}).
Consider the $\sl(1,n)$-modules\def\kop{K^\mathrm{op}}
\def\ksl{K^\mathrm{sl}}
\[
\ksl(\lambda) = \ind_{\sl_0\oplus\sl_{1}}^\sl L_\lambda, \qquad
\kop(\lambda) = \ind_{\sl_{-1}\oplus\sl_0}^\sl L_\lambda.
\]
There is a generic condition called \emph{typicality}
such that
$$\lambda\text{ is typical} \iff \ksl(\lambda) \text{ is irreducible}.$$
Moreover,\ $\lambda$ is typical iff $\ksl(\lambda)$ is projective
in the category of finite-dimensional $\gl(n)$-semisimple
$\sl(1,n)$-modules,
and the category of finite-dimensional $\gl$-semisimple
$\sl(1,n)$-modules with typical subquotients is semisimple.
For typical weights, we have $\kop_\lambda=\ksl(\lambda + 1\cdots1)$.
Note that if $\lambda\not\in\ZZ^n$, then the weight~$\lambda$ must be
typical.

Now consider certain $\W$-modules, the
\emph{big Kac modules}
$$
K'(\lambda) = \ind_{\W_{-1}\oplus\W_0}^\W L_\lambda = \ind_\sl^\W
\ind_{\g_{-1}\oplus\g_0}^\sl L_\lambda.
$$
We have $K'(\lambda) \isom \ind_\sl^\W \kop_\lambda = \ind_\sl^\W
\ksl(\lambda + 1\cdots1)$.
The $K'(\lambda)$ are indecomposable modules with highest weight
$\lambda$ with respect to~$\borel_2=\W_{-1}\oplus\borel_0$.

\section{Extensions}

In this section, we show the existence of certain non-split extensions
between Kac modules.  These will be realised as quotients of big Kac
modules.  More precisely, let $\M$ be the category of $\W(n)$-modules
that are direct sums of 
finite-dimensional semisimple $\gl(n)$-modules and are
$\sl(1,n)$-locally finite. 
\begin{claim}
(1)~$\M$ has enough projectives, i.e., every module in~$\M$ is a
   quotient of a projective module in~$\M$; (2)~if $P$ is a finite-dimensional
   projective $\sl$-module, then $\ind_\sl^\W P$ is projective
   in~$\M$.
\end{claim}
\begin{proof}
The main point is that induction takes projectives to projectives.
All modules in this proof are assumed to be direct sums of simple
finite-dimensional $\gl(n)$-modules.

The category of locally finite $\sl$-modules has enough projectives.
This may be checked as follows.
In the category of semisimple locally finite
$\gl(n)$-modules, every module is projective and injective.  
Therefore, if $Z$ is any such $\gl(n)$-module, then $\ind_\gl^\sl(Z)$
is projective in the category of locally finite $\sl$-modules.  In
particular,
if $M_0$
is a locally finite $\sl(1,n)$-module, then $M_0$ will be a quotient
of the standard projective $\ind_\gl^\sl(M_0)$, and therefore the category
of locally finite $\sl$-modules has enough projectives.

Finally, if
$M\in\Ob\M$, then there will exist a projective cover $P_0\to\res_\sl
M$, therefore we can get $M$ as a quotient
$\ind_\sl^\W(P_0)\to\ind_\sl^\W(M)\to M$, where the last arrow is the
canonical homomorphism $\ind_\sl^\W(M)\to M$, which is given by $u
\tensor m \mapsto um$.  This
proves~(1).

For~(2), we first note that, in this case, a locally finite
$\sl$-module will in fact be a direct sum of finite-dimensional
$\sl$-modules, and, therefore, if $P$ is projective in the category of
finite-dimensional $\sl$-modules, then $P$ is projective in this
slightly larger category.
Then induction makes
$\ind_\sl^\W(P)$ a projective $\W$-module.
\end{proof}
Because of the claim, since there are enough projectives, any module
in~$\M$ has a projective resolution.  If
\[
\cdots \to P_{-1} \to P_0 \to M \to 0
\]
is a projective resolution of a module~$M$, we can calculate
$\Ext_\M(M,N)$ as the cohomology of the complex
\[
\Hom_\M(P_0,N) \to \Hom_\M(P_{-1},N) \to \Hom_\M(P_{-2},N) \to \cdots,
\]
where the maps are obtained by composing a homomorphism with the
appropriate map in the projective resolution.
The following proposition is standard (e.g., see~\cite{Jer}):
\begin{prop}
If $\lambda\neq\mu$, then
$\Ext^1_{U(\W)}\bigl(S(\lambda),S(\mu)\bigr)=\Ext^1_\M\bigl(S(\lambda),S(\mu)\bigr)$.
\end{prop}
\begin{proof}
Let $E$ be an extension of $S(\lambda)$ by $S(\mu)$, where
$\lambda\neq\mu$; then we need to show that the centre of $\gl(n)$ acts
semisimply.

Let $\pi\colon  \W \to \End_\CC(E)$ be the action of~$\W$.  Let $Z$ be the Euler vector
field, and consider $\pi(Z)\in\End_\CC(E)$.  The
Jordan-Chevalley decomposition gives $\pi(Z)=X+Y$, where $X$
(resp.~$Y$) is semisimple (resp.\ nilpotent).  Furthermore, $\ad(X)$
and $\ad(Y)$ are respectively semisimple and nilpotent, they commute,
the decomposition is unique, and each of them maps $\pi(\W)$ into $\pi(\W)$.
Note that $\ad(Z)$ is semisimple, hence $\ad(\pi(Z))\in\End(\pi(\W))$ is
semisimple.  Then $\ad(Y)=\ad(\pi(Z))-\ad(X)$, the difference of two
commuting semisimple endomorphisms, so $\ad(Y)$ also
acts semisimply on $\pi(\W)$, therefore $\ad(Y)$ restricted to
$\pi(\W)$ is zero.

Therefore we see that
the nilpotent part~$Y$ of the action of the Euler vector
field is a nilpotent $\W$-module endomorphism of~$E$.  Since
$Y$ is nilpotent, we have $Y|_{S(\mu)}=0$ and
$Y\colon  S(\lambda) \to S(\mu)$.  If
$\lambda\neq\mu$, the only morphism $S(\lambda)\to S(\mu)$ is zero. 
\end{proof}

Suppose that $\lambda\notin\ZZ^n$.
Then $\lambda$ is typical for~$\sl(1,n)$, and
we have $K'(\lambda) \isom \ind_\sl^\W \kop_\lambda = \ind_\sl^\W
\ksl(\lambda + 1\cdots1)$.  
Such big Kac modules are therefore projective, according to the above
claim.
Now let $M$ be a finite-dimensional $\W(n)$-module all of whose simple
subquotients are Kac modules with nonintegral highest weight.  Then we
construct a projective resolution of~$M$ consisting of big Kac modules.
Construct the standard Koszul resolution
$$
\cdots \to P_{-2} \to P_{-1} \to P_0 \to \CC \to 0,
$$
where $P_{-i} = U(\W) \tensor_{U(\sl)} \bigwedge^i(\W/\sl)$, which is exact.
Tensoring by $M$, and using
$\ind_\p^\g(X)\tensor_\CC Y =
\ind_\p^\g(X\tensor_\CC\res Y)$,
we get a
projective resolution
$$
\cdots \to \ind_\sl^\W\bigl((\W/\sl)\tensor_\CC M) \to
\ind_\sl^\W M \to M \to 0
$$
from which we would like to calculate extensions of~$M$.
It is projective, since each $\bigwedge^i(\W/\sl)\tensor M$ is a
$\sl$-module all of whose simple subquotients are typical, hence
projective, Kac modules.

Note that the following results are vacuous unless $n\ge3$.
\begin{thm}\label{generic_exts}
Let $\lambda\in\Lambda^+\setminus\ZZ^n$, and let $\alpha$ be
a root of~$(\W/\sl)_1$ (i.e., $\alpha$ is a root of the complement of
$\sl$ in~$\W$ and $\sum_{i=1}^n\alpha_i=1$) such that $\lambda+\alpha\in\Lambda^+$.
Then $\Ext^1(K_\lambda,K_{\lambda+\alpha})\neq0$.
Moreover, the dimension of the space of extensions
is equal to
the multiplicity of $K_{\lambda+\alpha}$ in~$(\W/\sl)\tensor K_\lambda$.
\end{thm}

\begin{cor}\label{lots_of_exts}
Let $\lambda\in\Lambda^+$ be any weight, and let $\alpha$ be a
root of~$(\W/\sl)_1$.  Then
$\dim\Ext^1(K_\lambda,K_{\lambda+\alpha})\ge
[(\W/\sl)\tensor K_\lambda : K_{\lambda+\alpha}].$
\end{cor}
\begin{proof}[Proof of Theorem~\ref{generic_exts}]
The condition on~$\lambda$ ensures that all the Kac modules
involved are $\sl$-typical.
By the remarks above, in~$\M$ there is a projective resolution
$$\cdots \to P_{-1} \to P_0$$
of~$K_\lambda$, with
$P_{-i}=\ind_\sl^\W \bigl(\bigwedge^i(\W/\sl)\tensor K_\lambda\bigr)$.
The $\sl$-module $\bigwedge^i(\W/\sl)\tensor K_\lambda$ is
typical, therefore semisimple, so it is a direct sum of Kac modules.

We can calculate $\Ext_\M(K_\lambda,N)$ as the cohomology of the complex
$$\Hom_\M(P_0,N)\to\Hom_\M(P_{-1},N)\to\cdots.$$
Now,
$$\Hom_\W(\ind_\sl^\W K_\mu,K_\xi)=\Hom_\sl(K_\mu,K_\xi)=
\begin{cases}
0 & \text{if $\mu\neq\xi$},\\
\CC & \text{if $\mu=\xi$.}
\end{cases}$$
Therefore, if we substitute
$N=K_{\lambda+\alpha}$ into the above complex, we get
$$\Hom_\sl(K_\lambda,K_{\lambda+\alpha})\to\Hom_\sl(\W/\sl\tensor K_\lambda,K_{\lambda+\alpha})
\mathrel{\smash{\mathop{\to}\limits^\delta}}
\Hom_\sl(\textstyle\bigwedge^2(\W/\sl)\tensor K_\lambda,K_{\lambda+\alpha})\to\cdots,$$
so there are no coboundaries, and we claim that every
$f\in\Hom_\sl(\W/\sl\tensor K_\lambda,K_{\lambda+\alpha})$ gives a
cocycle.  Indeed,\ $(\delta f)(x)=f(dx)$, where
$d\colon  \bigwedge^2(\W/\sl)\tensor K_\lambda \to (\W/\sl)\tensor K_\lambda$,
and it is easy to verify that there are no nonzero $\sl$-module
maps $\bigwedge^2(\W/\sl)\tensor K_\lambda \to K_{\lambda+\alpha}$:
define the height of a weight~$(\lambda_1,\ldots,\lambda_n)$ to be
$\sum_{i=1}^n\lambda_i$.  Then the height of any $\nu$ such that
$\bigwedge^2(\W/\sl)\tensor K_\lambda = \bigoplus_\nu K_\nu$ is
greater than or equal to $\height{\lambda}+2$, while
$\height{\lambda+\alpha}=\height{\lambda}+\height{\alpha}=\height{\lambda}+1$.
This weight calculation shows that $K_{\lambda+\alpha}$ simply does not occur among
the $K_\nu$, and we are done.

Since all the Kac modules are typical, we can calculate the decomposition
of $(\W/\sl)\tensor \ksl_\lambda$ into a direct sum of Kac modules as the decomposition
of $(\W/\sl)\tensor L_\lambda$ into a direct sum of $\gl(n)$-modules.  In particular,
for~$\alpha$ as in Theorem~\ref{generic_exts}, the multiplicity
of~$K_{\lambda+\alpha}$ is nonzero.
\end{proof}

\begin{proof}[Proof of Corollary~\ref{lots_of_exts}]
Consider the cohomology
$H^1\bigl(W,\Hom(K_{\lambda+(t\cdots t)},K_{\lambda+\alpha+(t\cdots t)})\bigr)$
as $t\in\CC$ varies.  The complex computing this cohomology
is finite-dimensional, and shifting the weights by $t$ does not
change the dimension of the components.  We can therefore view it as a
complex with fixed terms with a differential that depends polynomially on~$t$.
By Theorem~\ref{generic_exts},
$\dim H^1=\bigl[(\W/\sl)\tensor K_\lambda : K_{\lambda+\alpha}\bigr]$ for
generic values of~$t$.  By semicontinuity, $\dim H^1$ can only increase
at special values.
\end{proof}

\section{Blocks and wildness}\label{sec:blocks_for_W}

A block of an Abelian category~$\M$ is defined to be an indecomposable
full Abelian subcategory that is a direct summand.
Given a subset~$\Gamma\subseteq\Irr\M$, we denote by $\M(\Gamma)$ the
full subcategory of~$\M$ consisting of objects all of whose simple subquotients
are in~$\Gamma$.

In this section, we take $\M$ to be the category of all finite-dimensional
representations of~$\W(n)$.  All objects of~$\M$ are generalised
weight modules, therefore there is a decomposition
according to weight and parity:
define parameters~$t=(\bar t,p(t))\in(\CC/\ZZ)\times\ZZ_2$.

We need the following facts about representations of Lie
superalgebras:
\begin{lemma}\label{weight_splitting}
Let\/ $\g$ be a Lie superalgebra (with a fixed Cartan\/ $\h$), and\/ $M$
be a generalised weight module.  Then\/ $M=\bigoplus_{t\in\h_0^*/Q}M(t)$ as a\/ $\g$-module,
where\/ $Q$ is the root lattice, and\/
$M(t)=\bigoplus\{\,M^{(\lambda)} \mid \lambda\in t\,\}\subseteq M$.
\end{lemma}
\begin{proof}
Consider a generalised weight module $M$, so that $M=\bigoplus_\alpha M^{(\alpha)}$ as a vector space
(recall that, for $\alpha\in\h_0^*$, the subspace $M^{(\alpha)}\subseteq M$
is simply $\{\,v\in M \mid \forall H\in\h_0\;\exists n\in\ZZ_+\;(H-\alpha(H))^nv=0\,\}$ in this case).
We have
the following simple fact:
if $\lambda,\,\mu\in\h_0^*$, then
$U(\g)^{(\lambda)}M^{(\mu)}\subseteq M^{(\lambda+\mu)}$.
An immediate
consequence is that the $M(t)$, defined above, are submodules.
\end{proof}

\begin{lemma}\label{parity_splitting}
Let\/ $\g$ be as in the previous lemma, and assume that,
for every\/ $\alpha\in\Delta$, we have\/ $\dim\g^{(\alpha)}=(0| k)$\/ or $(k|0)$.  Suppose\/ $M$ is a generalised weight module whose support is contained
in a single\/ $Q$-coset\/ $t=\lambda+Q$.  Then there exists
a parity function\/ $\sigma\colon  M\to M$, commuting with the
action of\/ $\g$, such that\/ $M=M'\oplus M''$, defined by\/
$M'=\{\,v\in M \mid \sigma(v)=v\,\}$, \ $M''=\{\,v\in M \mid \sigma(v)=-v\,\}$.
\end{lemma}
\begin{proof}
First, define a parity $p\colon  \Delta\to\ZZ_2$ by $p(\alpha)=p(X_\alpha)$
for some $X_\alpha\in\g^{(\alpha)}$; this is well-defined, by our hypothesis.
It extends linearly to a function $p\colon  Q\to\ZZ_2$.

Now, suppose $M$ is a generalised weight module whose support is contained
in a single $Q$-coset $t=\lambda+Q$.  Shift the parity function to
$p\colon  t \to \ZZ_2$ by setting
$p(\lambda+\alpha)=p(\lambda)+p(\alpha)$, where $p(\lambda)\in\ZZ_2$
is fixed arbitrarily.  Consider the linear map $\sigma\colon  M\to M$,
uniquely defined by requiring that, if $v\in M^{(\mu)}_d$, then
$\sigma(v)=(-1)^{p(\mu)-d}v$.  Finally, note that, if
$X_\alpha\in\g^{(\alpha)}$, then $X_\alpha v \in
M^{(\mu+\alpha)}_{d+p(\alpha)}$, so $\sigma(X_\alpha
v)=(-1)^{p(\mu)-d}X_\alpha v = X_\alpha\sigma(v)$.  Therefore $\sigma$
commutes with the action of $\g$, and $M$ breaks up into a direct sum
of two submodules, determined by this new parity.
\end{proof} 

\begin{prop}
Let $\g=\W(n)$, and let $M$ be a generalised weight $\g$-module.  Then 
$M$ decomposes, as a $\g$-module, into a direct sum
\[
M = \bigoplus_{\bar t\in\h_0^*/Q}M(\bar t) = \bigoplus_{\bar
t\in\h_0^*/Q}M(\bar t)'\oplus
M(\bar t)'',
\]
defined by Lemmas \ref{weight_splitting} and~\ref{parity_splitting}.
\end{prop}
Since we are interested in finite-dimensional modules, we see that
$\Gamma_{\bar t}$ consists of simple modules whose highest weight~$\lambda$
satisfies $\lambda_1 \equiv \bar t\pmod{\ZZ}$, where $\bar t\in\CC/\ZZ$.  The
categories $\M(\Gamma_{\bar t}')$ and~$\M(\Gamma_{\bar t}'')$ are equivalent via
parity-reversal, and for what follows we will not worry about whether
a highest weight is ``even'' or ``odd'', that being determined
according to Lemma~\ref{parity_splitting}.

\begin{thm}\label{block_decomposition}
The decomposition
$$\M = \bigoplus_{t\in(\CC/\ZZ)\times\ZZ_2}\M_t$$
is the block decomposition of~$\M$, i.e.,
the categories~$\M_t=\M(\Gamma_t)$ are indecomposable.
\end{thm}
\begin{proof}
We define a relation~$\myto$ on the set of highest weights, such that
$\lambda\myto\mu$ implies the existence of a finite-dimensional
indecomposable having both $S(\lambda)$ and $S(\mu)$ as subquotients,
and, consequently, that $S(\lambda)$ and $S(\mu)$ are in the same
$\M$-block.
Finally, we show that if
$S(\lambda),S(\mu)\in\Gamma_t$, i.e., if
$\lambda_1\equiv\mu_1\pmod{\ZZ}$, then we can get from $\lambda$
to~$\mu$ with a finite number of intermediate steps.

If $\bar t \not\in \ZZ$, then all Kac modules in~$\M_t$ are simple.
In that case, by Theorem~\ref{generic_exts}, there exists a
nonsplit extension of $S(\lambda)$ by $S(\lambda+\veps_i)$ as
long as $\lambda+\veps_i\in\Lambda^+$.

For the atypical case $\bar t\in\ZZ$, Corollary~\ref{lots_of_exts}
still ensures that there is an indecomposable module with subquotients
$K(\lambda)$ and $K(\lambda+\veps_i)$, and hence both $S(\lambda)$
and $S(\lambda+\veps_i)$ are subquotients, for~$1\le i\le n$ such
that $\lambda+\veps_i\in\Lambda^+$.

The relation $\lambda\myto\mu$ is now defined as follows:
set
$\lambda\myto\mu$ if $\mu=\lambda+\veps_i$ for some $1\le i\le n$.
We have established that $\lambda\myto\mu$ implies that $\lambda$ and
$\mu$ are in the same block.  Finally, it is clear that the closure
of~$\myto$ to an equivalence relation on~$\Lambda^+$ has equivalence
classes which are exactly the cosets $(\lambda+Q)\cap\Lambda^+$.
\end{proof}

\begin{thm}\label{block_wildness}
Each block~$\M_t$ is wild.
\end{thm}
\begin{proof}
Let $\lambda\in\Lambda^+$ be a weight such that
$K(\lambda)\in\Ob(\M_t)$,\ $K(\lambda)$ is simple,\ $\lambda+\alpha\in\Lambda^+$ for
every root $\alpha$ of $(\W/\sl)_1$, and all the $K(\lambda+\alpha)$
are also simple.  For example, any sufficiently dominant weight, i.e., $\lambda_1 \gg
\lambda_2 \gg \cdots \gg \lambda_n$, with $\lambda_i \equiv \bar
t\pmod{\ZZ}$, will do.  Then, by Corollary~\ref{lots_of_exts}, there
exists a nontrivial extension of $K(\lambda)=S(\lambda)$ by $K(\lambda+\alpha)=S(\lambda+\alpha)$.

Therefore, the $\Ext$-quiver of each block
contains a subquiver
consisting of a vertex~$\lambda$ with
arrows from it to~$\lambda+\alpha$ for each root $\alpha$
of~$(\W/\sl)_1$.  Not counting multiplicities, there are
$3\binom{n}{3}+n$ such roots, namely, $\veps_i+\veps_j-\veps_k$ with
$1\le i,\,j,\,k \le n$ and $i\neq j$.  Since $n\ge3$, we always have
$3\binom{n}{3}+n > 5$, and the resulting quiver is already wild (q.v. Proposition~\ref{wildness}):
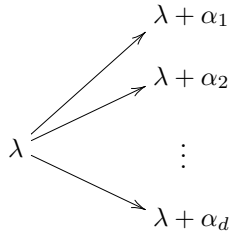
\[\xymatrix @R=9pt {
& **[r] {\lambda+\alpha_1} \\
& **[r] {\lambda+\alpha_2} \\
{\lambda} \ar[uur] \ar[ur] \ar[dr] & 
{\setbox0=\hbox{$\lambda+\alpha_d$}\kern\wd0\vdots\kern.5\wd0} \\
& **[r] {\lambda+\alpha_d}
}\]
\end{proof}

The original version of this paper used a much more cumbersome method
to establish the block decomposition.  The advantage of the present
approach is that much less calculation is necessary.
Also, it seems that the argument in
Theorem~\ref{lots_of_exts} may be refined to calculate \emph{all}
extensions between two (generic) Kac modules, and that the projective
resolution will give an easy proof of
\begin{conjecture}
All blocks not containing the trivial representation are equivalent.
\end{conjecture}

\section{The centre of $U\bigl(\W(n)\bigr)$}


Consider the standard embedding $\sl(1,n) \hookrightarrow \W(n)$.
Let $\h=\h_0\subset\gl(n)\subset\sl(1,n)$
be the standard Cartan, and fix a
triangular decomposition $\g = \n_- \oplus \h \oplus \n_+$ of~$\W(n)$
such that $\borel_+=\n_+\oplus\h$ is also the positive Borel of a
triangular decomposition of~$\sl(1,n)$.

Define the Verma modules    \def\msl{{^\sl\!M}}\def\lsl{{^\sl\!L}}
\begin{align*}
\msl(\lambda) &= U(\sl) \tensor_{U(\borel_+)} \CC_\lambda\\
M(\lambda) &= U(\W) \tensor_{U(\borel_+)} \CC_\lambda.
\end{align*}
Induction by stages gives $M(\lambda) = U(\W) \tensor_{U(\sl)} \msl(\lambda)$.
Let $\lsl(\lambda)$ be the unique irreducible quotient of~$\msl(\lambda)$.

The following result is proved in~\cite{Musson}, Section~3:
\begin{lemma}\label{musson}
Let $S\subseteq\h_0^*$.  Then
\[
\bigcap_{\mu\in S}\Ann_{U(\sl)} \lsl(\mu) = 0
\]
if and only if $S$ is Zariski dense in~$\h_0^*$.
\end{lemma}
\begin{cor}\label{gorelik}
For any Zariski dense subset~$S\subseteq\h_0^*$ one has
$$\bigcap_{\mu\in S}\Ann_{U(\W)} M(\mu) = 0.$$
\end{cor}
\begin{proof}
Suppose $S$ is Zariski dense.
First of all, $\Ann\msl(\mu) \subseteq \Ann\lsl(\mu)$ and therefore
$\bigcap_\mu \Ann_{U(\sl)}\msl(\mu) = 0$ by Lemma~\ref{musson}.  Next, applying the
Poincar\'e-Birkhoff-Witt Theorem, produce a basis~$\{X_i\}$ of $U(\W)$ 
over~$U(\sl)$, and write $u\in\Ann_{\U(\W)}M(\mu)$ as
$u=\sum_i X_i Y_i$, with the $Y_i\in U(\sl)$.  Applying $u$
to~$v\in1\tensor\msl(\mu)$ gives
$$\sum_iX_iY_iv=0,$$
and, since the $X_i$ are linearly independent, each $Y_iv=0$ and
therefore $Y_i\in\Ann_{U(\sl)}\msl(\mu)$.  This shows that
$\Ann_{U(\W)}M(\mu)\subseteq U(\W)\Ann_{U(\sl)}\msl(\mu)$ and that
\[
\bigcap_\mu \Ann_{U(\W)}M(\mu) \subseteq
U(\W)\bigcap_{\mu}\Ann_{U(\sl)} \msl(\mu) = 0.
\]
\end{proof}

\begin{prop}
We have\/ $Z\bigl(U(\W(0,n))\bigr)=\CC$.
\end{prop}

This is an immediate corollary of Corollary~\ref{gorelik}
and
\begin{lemma}
If $u\in Z$, then there exists a scalar~$c\in\CC$ such that
$u-c\in\bigcap_{\mu\in\Gamma}\Ann M(\mu)$, where $\Gamma\subset\h^*$
is Zariski dense.
\end{lemma}

\begin{proof}
For any $\mu\in\h^*$, we have $\End_\g M(\mu)=\CC$, therefore $c$ acts
by some scalar $c_\mu$ on~$M(\mu)$.  It therefore also acts by
multiplication by $c_\mu$ on the finite-dimensional quotient
$V_\borel(\mu)$ of~$M(\mu)$.  Furthermore, if $V_\borel(\mu)$ and
$V_\borel(\mu')$ are in the same block, then $c_\mu=c_{\mu'}$.
However, by Theorem~\ref{block_decomposition}, any block~$\Gamma$ is
already a Zariski dense subset of~$\h^*$.
\end{proof}


\bibliographystyle{plain}
\bibliography{main}

\end{document}